\documentclass[onefignum,onetabnum]{siamart190516}

\usepackage{lipsum}
\usepackage{amsfonts}
\usepackage{graphicx}
\usepackage{epstopdf}
\ifpdf
  \DeclareGraphicsExtensions{.eps,.pdf,.png,.jpg}
\else
  \DeclareGraphicsExtensions{.eps}
\fi

\newsiamremark{remark}{Remark}
\newsiamremark{hypothesis}{Hypothesis}
\crefname{hypothesis}{Hypothesis}{Hypotheses}
\newsiamthm{claim}{Claim}

\headers{A Linear View on Shape Optimization}{S.~Schmidt and V.~Schulz}

\title{A Linear View on Shape Optimization
}

\author{Stephan Schmidt\thanks{Humboldt University, Berlin 
  (\email{s.schmidt@hu-berlin.de}).}
\and Volker H.~Schulz\thanks{Trier University 
  (\email{volker.schulz@uni-trier.de}).}
}

\usepackage{amsopn}

\usepackage{algpseudocode}
\usepackage{algorithm}
\usepackage{mathrsfs}
\usepackage{bm} 

\def\R {{\mathbb R}}

\def\N {{\mathbb N}}

\newcommand{\vhsbf}[1]{{#1}} 

\def\bb{\vhsbf{b}}

\def\bf{{\vhsbf{f}}}

\def\bs{{\vhsbf{s}}}

\def\bu{{\vhsbf{u}}}
\def\bV{{\vhsbf{V}}}
\def\bW{{\vhsbf{W}}}

\def\bx{{\vhsbf{x}}}
\def\bY{\vhsbf{Y}}

\def\bq{{\vhsbf{q}}}

\def\bc{{\vhsbf{c}}}
\def\bu{{\vhsbf{u}}}

\newcommand{\vhsbm}[1]{{#1}} 

\def\blambda{{\vhsbm{\lambda}}}

\newcommand{\CG}{\operatorname{CG}}

\newcommand{\ACS}{2pt}
\newcommand{\mat}[4]{{\arraycolsep\ACS
\left#1\begin{array}{@{}*{#2}{c}@{}}#4\end{array}\right#3}}

\usepackage{tikz}
\usetikzlibrary{arrows,shapes,patterns,calc}
\tikzstyle{format} = [draw, thin,fill=white, draw=black]
\usepackage{pgfplots}
\usepgfplotslibrary{fillbetween}
\pgfplotsset{compat=1.6}
\pgfplotsset{soldot/.style={color=black,only marks,mark=*}} 
\pgfplotsset{holdot/.style={color=black,fill=white,only marks,mark=*}}

\def\norm#1{\|#1\|}
\def\bei#1{\vrule width 0.4pt height 14pt depth 9pt
           \lower 8pt \hbox{$ _{\hbox{} #1}$}\!\!\!}
\def\grad{{\nabla}}
\newcommand{\STS}[1]{\textcolor{blue}{#1}}

\newcommand{\Lag}{{\mathscr{L}}} 

\usepackage{tikz}

\def\d{\operatorname{d}}
\def\D{\operatorname{D}}

\renewcommand{\div}{\operatorname{div}}
\newcommand{\tr}{\operatorname{tr}}

\newcommand{\dualp}[2]{\left({{#1}},{{#2}}\right)}
\newcommand{\inner}[2]{\left\langle{{#1}},{{#2}}\right\rangle}
\newcommand{\dm}[2]{{\operatorname{d_M} {#1} {#2}}} 
\newcommand{\ddm}[2]{{\operatorname{d_M^2} {#1} {#2}}} 

\usepackage{calc}
\newcommand{\eqcolon}{\mathrel{\resizebox{\widthof{$\mathord{=}$}}{\height}{ $\!\!=\!\!\resizebox{1.2\width}{0.8\height}{\raisebox{0.23ex}{$\mathop{:}$}}\!\!$ }}}
\newcommand{\coloneq}{\mathrel{\resizebox{\widthof{$\mathord{=}$}}{\height}{ $\!\!\resizebox{1.2\width}{0.8\height}{\raisebox{0.23ex}{$\mathop{:}$}}\!\!=\!\!$ }}}

\ifpdf
\hypersetup{
  pdftitle={An Example Article},
  pdfauthor={D. Doe, P. T. Frank, and J. E. Smith}
}
\fi

\begin{document}

\maketitle

\begin{abstract}
Shapes do not define a linear space. This paper explores the linear structure of deformations as a representation of shapes. This transforms shape optimization to a variant of optimal control. The numerical challenges of this point of view are highlighted and a novel linear version of the second shape derivative is proposed leading to particular algorithms of shape Newton type.
\end{abstract}

\begin{keywords}
  Shape optimization, PDE constrained optimization, Numerical methods
\end{keywords}

\begin{AMS}
  49M15, 49M41, 49Q10
\end{AMS}

\section{Introduction}
Shape optimization is a research topic of high interest and is also used in numerous fields of application. Here, we mention only some fundamental references \cite{Sokolowski-1992, Delfour-Zolesio-2001, PironneauxX2001, Haslinger-2003} and some selected applications \cite{JOTA-2009, AIAA-2013, BANGTSSON20031533}. As mentioned, e.g., in \cite{Delfour-Zolesio-2001}, shapes do not define a linear space. For example, there is no straight forward way to define the sum of two shapes. But linear spaces are the realm of standard optimization techniques--may it be finite dimensional or infinite dimensional vector spaces. There are various ways to deal with this situation: enforcing finite dimensional parameterization of shapes like in CAD (computer aided design), which limits the variability of shapes; exploring the manifold structure of shapes like in \cite{VHS-shape-Riemann, Welker2021}, which is restricted to certain shape classes; method of mapping \cite{Murat-1976a, Jameson-1988, BS-CSE-2012} which maps the boundary of straight reference domain to the deforming boundary under consideration and which again has limitation in the variability of the reachable shapes. 

Here, we investigate the space of deformations of an initial shape. While the method of mapping typically limits the freedom of shape variability for the sake of uniqueness of optimal solutions, deformations of the whole computational domain enable full freedom in design. Another advantage is that the set of deformation is still a linear vector space. The downside of this approach is the lack of uniqueness of the optimal solutions, which poses challenges to numerical algorithms, remedies of which are also investigated in this paper.

This paper is organized as follows. \Cref{sec:shape deriv} recalls the classical definitions of first and second order shape derivatives.
Our main results are in
\cref{sec:main}, where we show that the classical shape derivatives possess a linear space interpretation. \Cref{mat-der-test} introduces a convenient way to derive necessary optimality conditions in weak form. A model problem is introduced in \cref{sec:model}, for which numerical results are given in the final \cref{sec:numres}.

\section{Shape Derivatives}\label{sec:shape deriv}
Here, we recall basic definitions of first and second shape derivatives. Let $d\in\N$ and $\tau>0$. We denote by $\Omega\subset \R^d$ a bounded domain with Lipschitz boundary $\Gamma:=\partial\Omega$ and by $J$ a real-valued functional depending on it.
Moreover, let $\{T_t\}_{t\in[0,\tau]}$ be a family of bijective mappings $T_t\colon \Omega\to\R^d$ such that $T_0=id$.
This family transforms the domain $\Omega$ into new perturbed domains $\Omega_t:=T_t(\Omega)=\{T_t(x)\colon x\in \Omega\}$ with $\Omega_0=\Omega$ and the boundary $\Gamma$ into new perturbed boundaries $\Gamma_t:=T_t(\Gamma)=\{T_t(x)\colon x\in \Gamma\}$ with $\Gamma_0=\Gamma$. If you consider the domain $\Omega$ as a collection of material particles, which are changing their position in the time-interval $[0,\tau]$, then the family $\{T_t\}_{t\in[0,\tau]}$ describes the motion of each particle, i.e., at the time $t\in [0,\tau]$ a material particle $x\in\Omega$ has the new position $x_t:=T_t(x)\in\Omega_t$ with $x_0=x$.
The motion of each such particle $x$ could be described by the \emph{velocity method}, i.e., as the flow $T_t(x):=\xi(t,x)$ determined by the initial value problem
\begin{equation}
\begin{split}
\frac{d\xi(t,x)}{dt}&=V(\xi(t,x))\\
\xi(0,x)&=x
\end{split}
\end{equation}
or by the \emph{perturbation of identity}, which is defined by $T_t(x):=x+tV(x)$ where $V$ denotes a sufficiently smooth vector field. We will use the perturbation of identity throughout the paper. The \emph{Eulerian derivative} of $J$ at $\Omega$ in direction $V$ is defined by
\begin{equation}
\label{eulerian}
dJ(\Omega)[V]:= \lim\limits_{t\to 0^+}\frac{J(\Omega_t)-J(\Omega)}{t}.
\end{equation}
The expression $dJ(\Omega)[V]$ is called the \emph{shape derivative} of $J$ at $\Omega$ in direction $V$ and $J$ \emph{shape differentiable} at $\Omega$ if for all directions $V$ the Eulerian derivative (\ref{eulerian}) exists and the mapping $V\mapsto dJ(\Omega)[V]$ is linear and continuous. For a thorough introduction into shape calculus, we refer to the monographs \cite{Delfour-Zolesio-2001,Sokolowski-1992}.
In particular, \cite{Sokolowski-1996} states that shape derivatives can always be expressed as boundary integrals due to the Hadamard structure theorem.
The shape derivative arises in two equivalent notational forms:
\begin{align}
\label{dom_form}
dJ_\Omega[V]&:=\int_\Omega F(x)V(x)\, dx  &\text{(domain formulation)}\\
\label{bound_form}
dJ_\Gamma[V]&:=\int_\Gamma f(s)V(s)^\top n(s)\, ds  &\text{(boundary formulation)}
\end{align}
where $F(x)$ is a (differential) operator acting linearly on the perturbation vector field $V$ and $f\colon\Gamma\to\R$ with
\begin{equation}
dJ_\Omega[V]=dJ(\Omega)[V]=dJ_\Gamma[V].
\end{equation}
The boundary formulation (\ref{bound_form}), $dJ_\Gamma[V]$, acting on the normal component of $V$ has led to the interpretation as tangential vector of a corresponding shape manifold in \cite{VHS-shape-Riemann}.

The second shape derivative is defined as the first shape derivative of a shape derivative. That necessitates a second flow movement after the first flow movement. If we denote a second flow induced by a vector field $W$ by $T^W_s$, then the second shape derivative is defined as
\[
d^2J[W,V]:=\lim\limits_{s\to 0^+}\frac1s\left(\lim\limits_{t\to 0^+}\frac{J(T^W_s(\Omega_t)-J(T^W_s(\Omega))}{t} - \lim\limits_{t\to 0^+}\frac{J(\Omega_t)-J(\Omega)}{t}\right)
\]
In \cite{Delfour-Zolesio-2001} concrete expressions for the second shape Hessian are discussed in the case of the velocity method.

In addition to these basic definitions, usually also concepts like local shape derivative and material derivative are needed in order simplify the shape calculus. We use that in the sections below and refer to \cite{Berggren-2010} for the specific usage in computations. 

\section{Shape Optimization Algorithms and Convergence in Deformation Vector Space}
\label{sec:main}
A standard shape optimization algorithm applies the same shape deformation $T_t$ used in the definition of the shape derivative in order to generate iterated shapes via
\[
\Omega^{k+1}=T_{t^k}^{V^k}(\Omega^k)\, ,
\]
where $V^k$ is related to the shape derivative $dJ(\Omega^k)[W]$. The earliest--to the knowledge of the authors--convergence result for this type of algorithm can be found in \cite{Hintermueller-2005}. There, all $V^k$ are assumed to be from a Hilbert space ${\cal H}$, the first and second order shape derivative is assumed to exist, the second shape derivative is assumed to be uniformly bounded and the step vector field $V^k$ is defined by
\begin{equation}\label{algo-Hintermueller-crit}
b_k(V^k,W)=-dJ(\Omega^k)[W]\, ,\ \forall W\in {\cal H}
\end{equation}
where $b_k(.,.)$ are assumed to be uniformly lower and upper bounded bilinear forms in ${\cal H}$. If the step length $t^k$ is chosen according to a back tracking line search, and $J$ is assumed bounded from below, then \cite[Theorem 3.1]{Hintermueller-2005} shows $\lim_{k\to\infty}J(\Omega^k)=J^*$ for some $J^*\in\R$ and $\lim_{k\to\infty}\norm{V^k}=0$. 

This is a very general convergence result. However, it does only state converge in the objective function and does not guarantee the existence of a limiting shape. Nor does it quantify the convergence speed. The key problem for further discussions is the proper notion of the distance in shape spaces. In \cite{VHS-shape-Riemann}, shapes are characterized as elements of a Riemannian manifold and the concepts of optimization on manifolds are carried over to shape manifolds. Within this framework, convergence analysis for steepest descent, Newton and Quasi-Newton methods is presented. The major drawback with this approach is, however, the requirement of $C^\infty$ smoothness of the shapes under investigation, which does not really fit into a typical finite element framework from standard PDE numerics.

Here, we consider only shapes, which result from deformations of an initial shape and identify the resulting shapes with the deformation mapping producing them. This is related to the concept of pre-shapes as discussed in detail in \cite{preshapea,preshapeb}, but we do not want to complicate the discussions unnecessarily with this concept. Starting from an initial shape $\Omega^0$, here we assume to work only with deformed shapes of the type $T(\Omega^0)=\{T(\bx) : \bx\in\Omega^0\}$, where $T:{\cal D}\to \R^d$ is a mapping defined on the hold-all-domain ${\cal D}$ and thus an element of a natural vector space. A central requirement for the mappings representing shapes is that the respective mapping $\bx\mapsto T(\bx)$  is invertible. Otherwise, severe problems like self intersecting discretization meshes or intersecting shape boundaries are to be expected. Furthermore, we want to concatenate shapes representing mappings, where we have to make sure that we only concatenate mappings ${\cal D}\to {\cal D}$. Thus, those $T$ representing shapes are a subset of all possible $T:{\cal D}\to \R^d$.

We choose the framework of perturbation of identity and consider shape optimization methods, which update a shape by a vector field $V$, which we interpret as a mapping $V:{\cal D}\to \R^d$. 
Thus, we rewrite a standard shape update rule for updating a shape $\Omega^k=T^k(\Omega^0)$ by a vector field $V:{\cal D}\to \R^d$ in terms of the respective deformation mapping as
\[
\Omega^{k+1}=\Omega^k+V(\Omega^k) \Leftrightarrow
T^{k+1}(\Omega^0) = T^{k}(\Omega^0)+V\circ T^{k}(\Omega^0)
=(I+V)\circ T^{k}(\Omega^0)\, ,
\]
where we observe that $I+V$ is invertible, if $\|V\|_{W^{1,\infty}}<1$ \cite{Allaire-2007}, and that  
$I+V:{\cal D}\to {\cal D}$, if $V$ is zero close to at the boundary of ${\cal D}$. This is implicitly checked in each shape optimization iteration. 
Therefore, a standard shape optimization algorithm written as an algorithm in deformation vector fields, may be written in the form of \cref{algo-shape-steepest}.
\begin{algorithm}
\caption{Steepest Descent Shape Deformation Optimization}
\label{algo-shape-steepest}
\begin{algorithmic}[0] 
\State $k:=0$, initialize $T^0=I,\varepsilon$
\Repeat 
\State determine $V^k$ solving $b_k(V^k,Z)=-dJ(\Omega^k)[Z]\, ,\ \forall Z\in {\cal H}$
\State Line search: find $t^k\approx\arg\min_tJ((I+tV^k)\circ T^k(\Omega^0))$
\State $T^{k+1}:=(I+t^kV^k)\circ T^k$
\State $k:= k+1$
\Until{$\norm{V^k}\le\varepsilon$}
\end{algorithmic}
\end{algorithm}

At the first glance, algorithm \ref{algo-shape-steepest} differs significantly from standard steepest descent optimization algorithm in vector spaces as discussed, e.g.~in \cite{Nocedal-Wright}, in view of the update rule $T^{k+1}:=(I+t^kV^k)\circ T^k$. But because of the special definition of the shape derivative, this is indeed a standard steepest descent method in deformation vector fields, as is shown in the following \Cref{algo_grad_deform}. First we introduce the notation. We interpret the shape objective as a mapping
$T\mapsto f(T):=J(T(\Omega^0))$ and rewrite the shape derivative in terms of the perturbation vector field $Z$. 
\begin{equation}\label{eq:sdasfd}
dJ(\Omega)[Z]=\frac{d}{dt}\bei{t=0+}J((I+tZ)(\Omega))
= \frac{d}{dt}\bei{t=0+}f((I+tZ)\circ T)
=f'(T)[Z\circ T]
\end{equation}
Equation \cref{eq:sdasfd} is the central observation initiating this paper. It means that the standard shape derivative can be seen as a usual directional derivative for shape deformations, where the local directions $Z$ are pulled back to $Z\circ T$.
We would like to find a gradient related vector, which we call $\nabla J$ satisfying relation (\ref{algo-Hintermueller-crit}) for $Z$ reformulated in this notation as
\begin{equation}\label{algo_deT_gradJ}
b(\nabla J(T), Z)_T=dJ(\Omega)[Z], \quad \forall Z \in {\cal H}\, . 
\end{equation}  
The bilinear form $b(.,.)_T$ may be chosen differently at each point $T$ in the deformation vector space, which is indicated by the index $T$. We observe that $\nabla J$ cannot be interpreted as a gradient of the function $f$ defined on the Hilbert space of deformations $({\cal H},g)$ with scalar product $g$ in the following way
\begin{equation}\label{gradf-def}
g(\nabla f(T),Z)=f'(T)Z\, , \forall Z \in {\cal H}\, .
\end{equation}
However, there holds a revealing relation between both as formulated in \Cref{algo_grad_deform}.

\begin{theorem}\label{algo_grad_deform}
Let $g$ be a scalar product on the linear space of mappings ${\cal H}:=\{T \ \vert \ T:{\cal D}\to\R^d\}$  such that $({\cal H}, g)$ is a Hilbert space. 
Let $J:\Omega\to\R$, where $\Omega=T(\Omega^0)$, be a shape functional and $f:T\mapsto J(T(\Omega^0))$ be a related objective function in terms of a deformation $T\in{\cal H}$. Furthermore, we define $\nabla J$ by equation (\ref{algo_deT_gradJ}) and assume that all local bilinear forms $b(.,.)_T$  are related to an overall scalar product $g(.,.)$ on the space of deformations in the following natural way for any invertible $T\in{\cal H}$:
\begin{equation}\label{algo_scalar}
b(Z_1,Z_2)_T:= g(Z_1\circ T^{-1}, Z_2\circ T^{-1})\, .
\end{equation}
Then, there holds for any invertible $T\in{\cal H}$ the following relation between $\nabla J$ and the gradient of $\nabla f$, where the gradient $\nabla f$ is defined by the scalar product $g$:
\begin{equation}\label{algo_grad_grad}
\nabla f(T) = \nabla J(T)\circ T\,. 
\end{equation}
\end{theorem}
\begin{proof}
We observe from (\ref{algo_deT_gradJ}, \ref{algo_scalar})
\[
g(\nabla J(T)\circ T, Z\circ T)=f'(T)[Z\circ T],\  \forall Z\in{\cal H}\, .
\] 
Since $T$ is assumed invertible, this is equivalent to
\[
g(\nabla J(T)\circ T, Z)=f'(T)[Z],\  \forall Z\in{\cal H}\, ,
\] 
which defines $\nabla f$ according to (\ref{gradf-def}).
\end{proof}

Therefore, algorithm (\ref{algo-shape-steepest}) can be considered a steepest descent algorithm for the objective $f(X):=J(X(\Omega^0))$--potentially with variable metric in the case of remeshing or when using Newton's method.
Thus, all standard convergence considerations  apply here in the case of a finite dimensional Hilbert space $\cal H$ and natural generalizations as known in the optimal control community apply in the case of $\cal H$ being an appropriate infinite dimensional function space. Furthermore, we note  that condition (\ref{algo_scalar}) is quite natural, if the scalar product $b(.,.)_T$ is defined on the finite element mesh deformed by the deformation $T$.
%

In similar fashion, we can formulate a Taylor series in following theorem. 

\begin{theorem}\label{algo_shape_Taylor}
Let the assumptions of theorem \ref{algo_grad_deform} hold and assume for the vector field $V$ defined on ${\cal D}$ that $I+V:{\cal D}\to {\cal D}$ is invertible and that $(I+V)\circ T\in {\cal H}$ for a given $T\in {\cal H}$. Furthermore, we assume that the function $f$ defined in theorem \ref{algo_grad_deform} is three times differentiable with bounded third derivative. Then, there holds the following Taylor series
\begin{equation}\label{algo_Taylor_eq}
J(\Omega + V(\Omega))=J(\Omega)+b(\nabla J(\Omega),V)_T+\frac12 J^{\prime\prime}(\Omega)[V,V]+{\cal O}(\norm{V}^3)
\end{equation}
where ${\cal O}$ denotes the Landau symbol, the norm $\norm{.}$ is derived from the scalar product $b(.,.)_T$, and $J^{\prime\prime}(\Omega)[V,W]$ is the linear second shape derivative defined by
\begin{equation}\label{algo_linear_Hessian}
J^{\prime\prime}(\Omega)[V,W] := \frac{d}{ds_1}\bei{s_1=0+}\ \frac{d}{ds_2}\bei{s_2=0+}\ 
J((I+s_1V+ s_2W)(\Omega)).
\end{equation}
\end{theorem}
\begin{proof}
The assumptions guarantee the Taylor series for $f$, which is now rephrased. There holds
\begin{align*}
J(\Omega + &V(\Omega))=f(T+V\circ T) \\
&= f(T)+ f^\prime(T)[V\circ T]+\frac12 f^{\prime\prime}(T)[V\circ T,V\circ T]+{\cal O}(\norm{V\circ T}^3)\\
&= f(T)+ g(\nabla f(T),V\circ T)\\
&+\frac12  \frac{d}{ds_1}\bei{s_1=0+}\ \frac{d}{ds_2}\bei{s_2=0+} f(T+s_1V\circ T +s_2V\circ T)+{\cal O}(\norm{V\circ T}^3)\\
&=J(\Omega)+b(\nabla J(\Omega),V)_T+\frac12 J^{\prime\prime}(\Omega)[V,V]+{\cal O}(\norm{V\circ T}^3).
\end{align*}
The norm in the remainder term is derived from  the scalar product $g$. If we translate it to $b(.,.)_T$, we confirm the assertion.
\end{proof}


The Taylor series allows the construction of linear second order optimization methods, see below, and the investigation of well-posedness of shape optimization problems, although well-posedness considerations are simpler carried out on the shape boundary, either in terms of local parametrization as in \cite{EH-2012} or in terms of shape manifolds as in \cite{VHS-shape-Riemann}. It should be noted that the linear second shape derivative defined in theorem \ref{algo_shape_Taylor}, which is needed for the Taylor series, is symmetric and differs from the classical shape Hessian, which we recall here as
\[
d^2J(\Omega)[V,W]=\frac{d}{ds_1}\bei{s_1=0+}\ \frac{d}{ds_2}\bei{s_2=0+}\ 
J((I+s_1V)\circ(I+s_2W)(\Omega))
\]
The difference contains several nonsymmetric deritative terms of  ${\bV}\circ{\bW}$  arising in the classical shape Hessian. It is shown in \cite[section 6.5]{Delfour-Zolesio-2001} that the classical shape Hessian is symmetric, if ${\bV}\circ{\bW}={\bW}\circ{\bV}$ . The introduction of the linear second shape derivative above makes this somewhat artificial assumption no longer necessary.

\def\cH{{\cal H}}
\def\cN{{\cal N}}
\def\cR{{\cal R}}
\def\cL{{\cal L}}
\newcommand{\spy}{\,|\,} 

Based on the observations in \Cref{algo_shape_Taylor}, we formulate a shape Newton method
\begin{align*}
\Omega^{k+1}&=\Omega^{k}+V^k(\Omega^{k})\text{, where } V^k \text{ is defined by}\\
&J^{\prime\prime}(\Omega^k)[V^k,W]=-dJ(\Omega^k)[W]\, , \ \forall W\in{\cal H}
\end{align*}
which is equivalent to the formulation in deformation vector space
\begin{align*}
T^{k+1}&=T^{k}+V^k\circ T^{k}\text{, where } V^k \text{ is defined by}\\
&f^{\prime\prime}(T^k)[V^k\circ T^{k},W]=-f^{\prime}(T^k)[W]\, , \ \forall W\in{\cal H}
\end{align*}
with the assumption $\Omega^{k+1}=T^{k+1}(\Omega)$ and $\Omega^{k}=T^{k}(\Omega)$. Unfortunately, $J^{\prime\prime}$ and thus $f^{\prime\prime}$ have a nontrivial and even huge null-space, which cannot simply be ignored. If this would not be the case, then standard arguments would lead to locally quadratic convergence of the shape Newton method. In order to mitigate this problem, we use the Moore-Penrose pseudoinverse instead, which allows for locally quadratic convergence again.   
Since the convergence theory is more easily formulated for operators than for bilinear forms, we use again Riesz representation theorem in the Hilbert space $({\cal H},g)$ and define the gradient as in equation \cref{algo_grad_grad}
and with that also the Hessian operator $H(T)$ is defined by the operator representation of  $f^{\prime\prime}(T)$, i.e.
\[
g(H(T)V,W):=f^{\prime\prime}(T)[V,W]\, ,\ \forall V,W\in {\cal H}
\]
In this notation, the following iteration is  performed:
\begin{equation}\label{GN-it}
T^{k+1}=T^k + V^k\circ T^{k}\, ,\mbox{ where } V^k=-H(T^k)^+\grad f(T^k)
\end{equation}
Here $H(T^k)^+$ denotes the Moore-Penrose-pseudoinverse of $H(T^k)$. Specific aspects of the implementation are discussed in \Cref{algo-general-SVD} below. Iteration (\ref{GN-it}) can be shown to be locally quadratically convergent by applying well-known results from linear spaces as in \cite{DH-1979}.

The Moore-Penrose-pseudoinverse based on the SVD with the standard Euclidean scalar product is computationally inefficient. Therefore, we reformulate the Moore-Penrose-pseudoinverse as a least squares problem in the scalar product related to the Hilbert space we are working in and investigate a perturbation approach. We are considering the linear equation $H\bV=\bb$, where $H=H(T^k)$ is the Hessian operator and $\bb= \grad f(T^k)$ is the gradient at step $k$ of iteration (\ref{GN-it}), i.e., a Newton method based on the pseudoinverse.

\begin{theorem}\label{algo-general-SVD}
Let (${\cal H}, g$) be a Hilbert space with inner product $g$. We assume that the linear operator $A$ defined on ${\cal H}$ has closed range and is not necessarily invertible. When solving the equation $H\bV=\bb$ with $\bb\in\cR(H)$, we obtain $\hat{\bV}:=H^+\bb$ as the minimum norm solution, where $H^+$ is the Moore-Penrose pseudoinverse operator. Then, the vector $\hat{\bV}$ is also the unique solution of the optimization problem 
\begin{align}\label{algo-MPopt.1}
\min_\bV\ g(\bV,&\bV)\\\label{algo-MPopt.2}
\mbox{s.t. }\  H\bV&=\bb
\end{align}
Furthermore, if $H$ is self-adjoint in the scalar product $g$ and positive semidefinite, then the vector $\hat{\bV}$ can be computed as the limit of the solutions $\bV_\varepsilon$ of the  following family of linear-quadratic problems parameterized  by $\varepsilon >0$. For
\begin{equation}\label{algo-pseudo-eps}
\bV_\varepsilon :=\arg\min_\bV \frac12g(H\bV,\bV)-g(\bb,\bV) + \frac{\varepsilon}2g(\bV,\bV)
\end{equation}
there holds $\lim_{\varepsilon\to 0}\bV_\varepsilon = \hat{\bV}$.
\end{theorem} 

\begin{proof}
Instead of the frequently used definition of the pseudoinverse by an SVD, we use the equivalent variational definition  from
\cite[Definition (V) on page 45]{GroetschCxW1977}, which defines it as the minimum norm solution of equation (\ref{algo-MPopt.2}) in least squares reformulation. Thus any solution of the equation $H\bV=\bb$ can be written as
$\bV=H^+\bb+\bV_N=\hat{\bV}+\bV_N$, where $\bV_N\in\cN(H)$. We observe that $\hat{\bV}\in \cR(H^+)=\cR(H^*)$ \cite[Theorem 2.1.2]{GroetschCxW1977}, where $H^*$ is the adjoint operator. Therefore $\hat{\bV}\perp \bV_N$ for any $\bV_N\in\cN(A)$ \cite[Theorem 1.2.1]{GroetschCxW1977}, where $\perp$ means orthogonality in the $g$ scalar product. Using this in the objective (\ref{algo-MPopt.1}) gives the assertion of the first part.

The necessary conditions of problem (\ref{algo-pseudo-eps}) can be written as
\begin{align*}
g(H\bV_\varepsilon,\eta)-g(\bb,\eta)+\varepsilon g(\bV_\varepsilon,\eta) &= 0\, ,\ \forall \eta\in {\cal H}\\
\Leftrightarrow\quad H\bV_\varepsilon - \bb + \varepsilon \bV_\varepsilon &= 0
\end{align*}
We use the Ansatz $\bV_\varepsilon =: \hat{\bV}+\bY_\varepsilon$ for some $\bY_\varepsilon\in {\cal H}$. Since $H\hat{\bV}-\bb =0 $, we conclude
\begin{equation}\label{algo-ps-null}
(H+\varepsilon I)\bY_\varepsilon = -\varepsilon \hat{\bV}
\end{equation}
Since $\hat{\bV}\in \cR(H)$, because $H$ is self-adjoint, also $\bY_\varepsilon\in \cR(H)$ as a first consequence of (\ref{algo-ps-null}), for all $\varepsilon >0$. On the other hand, we observe in the same equation $H\bY_\varepsilon\to 0$ for $\varepsilon\to 0$ and thus $\lim_{\varepsilon\to 0 }\bY_\varepsilon\in \cN(H)$. Therefore
\[
\bV_\varepsilon - \hat{\bV} = \bY_\varepsilon \to 0\, ,\mbox{  for }\varepsilon\to 0\, .
\]

\end{proof}

In the context of a shape Newton method $H$ is the Hessian operator related to the linearized second shape derivative and $\bb$ is the gradient of the shape derivative.
\begin{corollary}\label{remark-Tichonov-KKT}
Due to the definition of the Hessian and the gradient operator, equation (\ref{algo-pseudo-eps}) is obviously equivalent to the formulation
\begin{equation}\label{algo-var-eps}
\min_\bV \frac12f^{\prime\prime}(T^k)[\bV\circ T^{k},\bV\circ T^{k}] - f^\prime(T^k)[\bV\circ T^{k}]+ \frac{\varepsilon}2g(\bV\circ T^{k},\bV\circ T^{k})
\end{equation}
which again, due to the definition of $f$, can be equivalently rephrased in terms of the shape functional as
\begin{equation}\label{algo-var-eps.s}
\min_\bV \frac12J^{\prime\prime}(\Omega^k)[\bV,\bV] - dJ(\Omega^k)[\bV]+ \frac{\varepsilon}2b_{\Omega_k}(\bV,\bV),
\end{equation}
which has to be solved in each shape Newton iteration. 
\end{corollary}

If $\varepsilon$ is reduced during the Newton iterations in the fashion $\varepsilon_k={\cal O}(\| \bV^k\|)$, then locally quadratic convergence of the resulting method can be expected from respective investigations in full Newton methods \cite{Deufl-Newton}. However, due to discretization effects, the  condition in theorem \ref{algo-general-SVD} that the gradient has to lie in the range of the Hessian may not be fulfilled exactly in the discretized equations. Therefore $\varepsilon$ may not be chosen arbitrarily small in practice.  
Formulations (\ref{algo-pseudo-eps}, \ref{algo-var-eps}, \ref{algo-var-eps.s}) are related to Tikhonov regularization. Thus, a similar result can be found in \cite[Corollary 2.3.8]{GroetschCxW1977} for the assumption that $H^*H$ is invertible, which does not hold here in general.

A different approach for the computation of the operator-vector product $H^+\bb$ can be performed by the usage of Krylov subspace methods, see, e.g., \cite{Hayami-2018, Brown-Walker-1997, Choi-Page-Saunders-2011}. However, the relevant publications have to be rephrased to take into account a general scalar product $g$ rather than the Euclidean scalar product. 

\section{Material Derivatives as Test Vectors}\label{mat-der-test}
So far, we have discussed derivatives of mappings from shapes to an objective criterion. Often, this mapping involves system equations in the form of partial differential equations, where it is reasonable to formulate the PDE explicitly as a constraint. 
Because of the intricacies of the shape calculus, it seems difficult to derive necessary optimality conditions for system model based shape optimization problems. Here, we introduce a simple step-by-step procedure for the derivation of the weak form of necessary conditions. We start out with the general necessary conditions of optimality for model based based optimization problems in the form
\begin{align*}
\min &J(\bu,\bq)\\
\mbox{s.t.}\ \dualp{\bc(\bu,\bq)}{\blambda}&=0\, , \ \forall \blambda\in \Lambda
\end{align*}
There, we build the Lagrangian 
\[
\Lag(\bu,\bq,\blambda):= J(\bu,\bq)+\dualp{\bc(\bu,\bq)}{\blambda}.
\]
and observe as necessary conditions of optimality
\begin{align}\label{turning_weak_short_1}
\frac{\partial}{\partial \bu}\Lag(\bu,\bq,\blambda)[\tilde{\bu}]=&\;0\, ,\ \forall \tilde{\bu} \in U\\\label{turning_weak_short_2}
\frac{\partial}{\partial \bq}\Lag(\bu,\bq,\blambda)[\tilde{\bq}]=&\;0\, ,\ \forall \tilde{\bq} \in Q\\\label{turning_weak_short_3}
\frac{\partial}{\partial \blambda}\Lag(\bu,\bq,\blambda)[\tilde{\blambda}]=&\;0\, ,\ \forall \tilde{\blambda} \in \Lambda,
\end{align}
where often $U=\Lambda$. Conditions (\ref{turning_weak_short_1} - \ref{turning_weak_short_3}) can be equivalently and jointly written as
\begin{equation}\label{turning_compact_KKT}
\frac{\partial}{\partial \bu}\Lag(\bu,\bq,\blambda)[\tilde{\bu}] +
\frac{\partial}{\partial \bq}\Lag(\bu,\bq,\blambda)[\tilde{\bq}] +
\frac{\partial}{\partial \blambda}\Lag(\bu,\bq,\blambda)[\tilde{\blambda}] =0
\, ,\ \forall (\tilde{\bu},\tilde{\bq},\tilde{\blambda}) \in U\times Q\times \Lambda
\end{equation}
since each individual equation from (\ref{turning_weak_short_1} - \ref{turning_weak_short_3}) follow from this equation by setting two out of the three vectors $\tilde{\bu},\tilde{\bq},\tilde{\blambda}$ to zero. We arrive at the same expression on the left hand side by explicitly denoting at the solution of the optimization problem that the state and adjoint depend on the optimization variable, i.e. $\bu=\bu(\bq), \blambda=\blambda(\bq)$ and calculating the total derivative of the Lagrangian with respect to $\bq$ as
\begin{align*}
&d\Lag(\bu(\bq),\bq, \blambda(\bq))[\tilde{\bq}]=\\
&\frac{\partial}{\partial \bu}\Lag(\bu,\bq,\blambda)[d\bu[\tilde{\bq}]] +
\frac{\partial}{\partial \bq}\Lag(\bu,\bq,\blambda)[\tilde{\bq}] +
\frac{\partial}{\partial \blambda}\Lag(\bu,\bq,\blambda)[d\blambda[\tilde{\bq}]].
\end{align*}
In this way, $d\Lag[\tilde{\bu}]$ depends additionally on the vectors $d\bu[\tilde{\bq}],d\blambda[\tilde{\bq}]$, which thus can be seen as somewhat fancy notations for the vectors $\tilde{\bu},\tilde{\blambda}$ in (\ref{turning_compact_KKT}). 

For standard model based opimization, this point of view is redundant. However, in the context of shape calculus, this approach gives a convenient guideline for deriving necessary conditions in the following steps
\begin{enumerate}
\item build a Lagrangian $\Lag(u,\Omega,\lambda)$ of the PDE constrained shape optimization problem based on state variable $u$ and adjoint $\lambda$
\item derive an expression for its shape derivative, where material derivatives \newline $\dm{u}{[\bV]},\dm{\lambda}{[\bV]}$ explicitly appear and denote this expression as \newline $d\Lag(u,\Omega,\lambda)[\dm{u}{[\bV]},\bV,\dm{\lambda}{[\bV]}]$ 
\item state necessary conditions in weak form as
\[
d\Lag(u,\Omega,\lambda)[\dm{u}{[\bV]},\bV,\dm{\lambda}{[\bV]}] = 0\, ,\ \forall (\dm{u}{[\bV]},\bV,\dm{\lambda}{[\bV]})
\in U\times H\times V
\]
where $U$ is the space for the state variables, $H$ is the Hilbert space for the admissible deformation vector fields of the shape, and $V$ is the space for the adjoint variables in the respective shape optimization problem.
\end{enumerate}
This process is carried out and exemplified in detail in the next section. Of course, in principle, also expressions for $d\Lag$ can be derived, where local shape derivatives instead of material derivatives appear. However, as noted in \cite{Berggren-2010}, only material derivatives can be considered from the same vector space as the state or adjoint variable. Furthermore, the existence of the material derivative of the adjoint might be in question. In this case, the derivation following steps 1-3 above can be considered only formal. Nevertheless, the resulting necessary conditions are in-line with a derivation avoiding material derivatives at all as discussed in \cite{Sturm-2015b}.

\section{Model Problem}\label{sec:model}

\begin{figure}[h!]
\begin{center}
\begin{tikzpicture}[]
  \draw [draw=black, fill=black!10] (0.0,0.0) rectangle (0.4\textwidth, 0.4\textwidth);
  \draw[draw=black, fill=black!30] (0.5*0.4\textwidth, 0.5*0.4\textwidth) ellipse ({0.25*0.4\textwidth} and {0.125*0.4\textwidth});
  \node[] at ({0.5*0.4\textwidth},{0.02\textwidth}) {$\Gamma_0$};
  \node[] at ({0.02\textwidth},{0.5*0.4\textwidth}) {$\Gamma_1$};
  \node[] at ({0.5*0.4\textwidth},{0.38\textwidth}) {$\Gamma_2$};
  \node[] at ({0.38\textwidth},{0.5*0.4\textwidth}) {$\Gamma_3$};
  \node[] at ({0.5*0.4\textwidth},{0.5*0.4\textwidth}) {$\Omega_0$};
  \node[] at ({0.85*0.4\textwidth},{0.85*0.4\textwidth}) {$\Omega_1$};
\end{tikzpicture}
\hspace{.05\textwidth}
\includegraphics[width=0.4\textwidth]{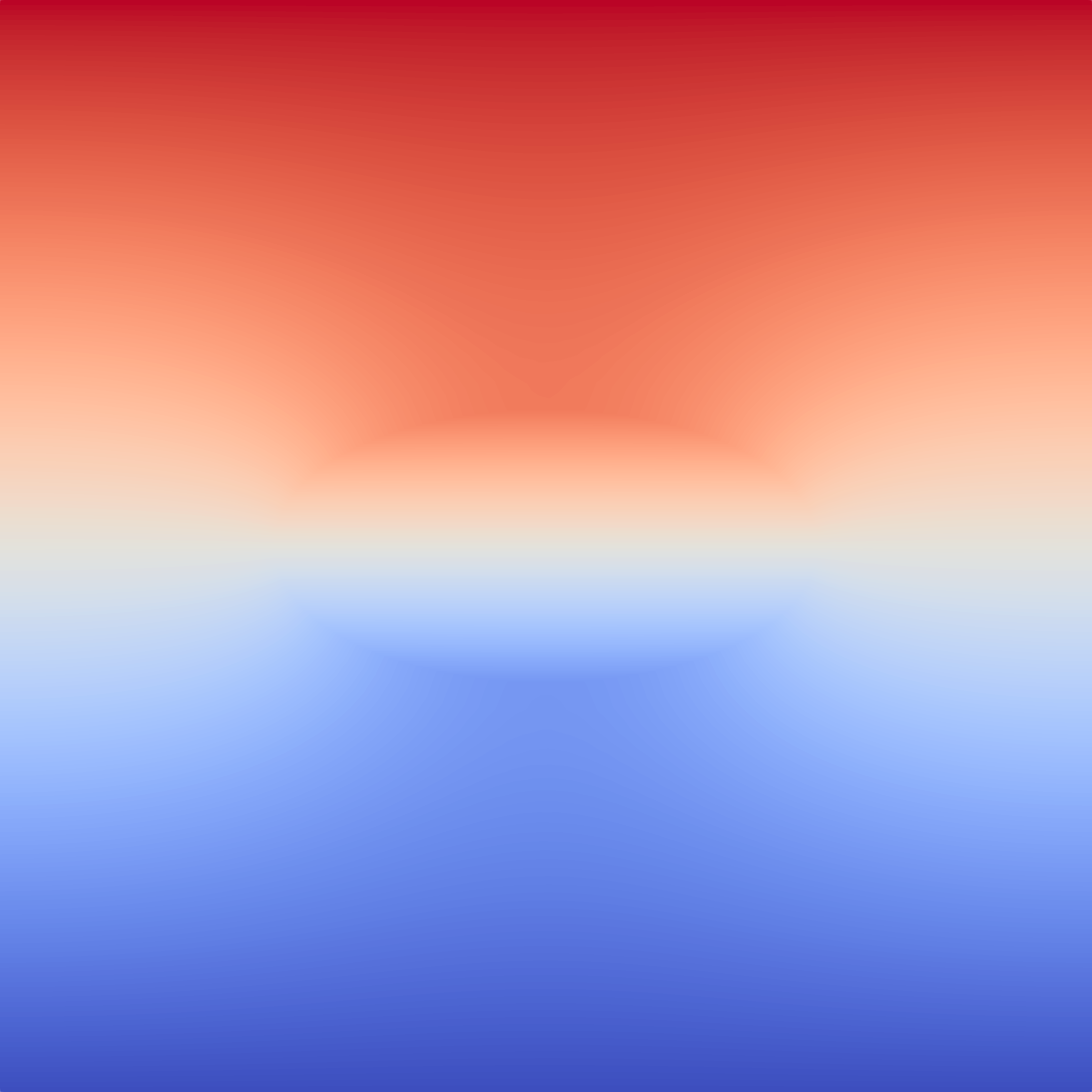}
\end{center}
\caption{The interface identification problem. Domain $\Omega$ to be reconstructed on the left, corresponding target state $z$ on the right.}\label{fig:InterfaceNewtonExample}
\end{figure}
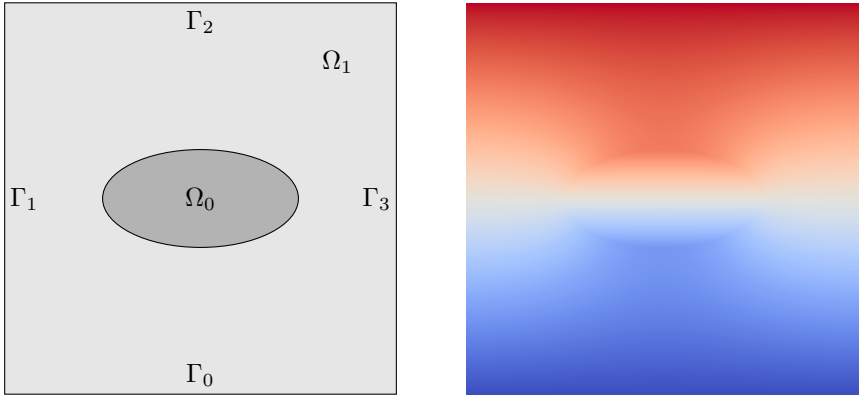

We consider the following model problem. We seek to reconstruct the shape of an inclusion $\Omega_0$ within the surrounding domain $\Omega_1$ with $[0,1]^2 \eqcolon \Omega = \Omega_0 \cup \Omega_1$ disjointly. The actual reconstruction is inspired by electrical impedance tomography and given by the following problem
\begin{equation}\label{eq:ShapeMotherProblem}
\begin{aligned}
  &\min_{(u, \Omega)}\quad \frac{1}{2} \int\limits_\Omega (u - z)^2 \d x + \frac{\alpha}{2} R(\Omega)\\
 \text{s.t.} & \text{ find } u \in u_0 + H_{(\Gamma_0\cup\Gamma_2)}^1(\Omega) \text{ such that} \\
 &\int\limits_{\Omega} \mu\inner{\nabla u}{\nabla v}_2 - f v \ \d x = 0 \quad \forall v \in H_{(\Gamma_0\cup\Gamma_2)}^1(\Omega),
\end{aligned}
\end{equation}
where $R(\Omega)$ stands for a regularization term, usually the perimeter of $\partial\Omega_0$ or the volume of $\Omega_0$. In the numerical examples below, we choose the volume and also $\alpha=10^{-6}$.
The affine Sobolev space $u_0 + H_{(\Gamma_0\cup\Gamma_2)}^1(\Omega)$, where $u_0(x_1, x_2) \coloneq x_2$ and $H_{(\Gamma_0\cup\Gamma_2)}^1(\Omega)=\{\varphi\in H^1(\Omega)\, :\, \varphi|_{\Gamma_0\cup\Gamma_2}=0\}$, incorporates the inhomogenous Dirichlet boundary conditions
\begin{align*}
u &= 0 \text{ on } \Gamma_0\\
u &= 1 \text{ on } \Gamma_2,
\end{align*}
while the ``do-nothing-condition'' naturally creates a homogenous Neumann condition on the remaining boundaries. 
These boundary conditions can be interpreted as putting an electric potential between the upper and lower side of the domain. The actual electric field lines would then be given by $\nabla u$. As such, the Neumann boundary conditions can be interpreted as no electric field lines entering or leaving the domain. The conductivity $\mu$ directly depends on the subdomains $\Omega_0$ and $\Omega_1$ via
\begin{align}\label{eq:InterfaceTwoStates}
  \mu(x) = \begin{cases}10^{-6}&, x\in \Omega_0\\1&, x\in\Omega_1\end{cases}.
\end{align}
We generate the desired potential $z$ by solving the problem beforehand for an elliptic inclusion $\Omega_0$ as shown in~\Cref{fig:InterfaceNewtonExample}.
Knowledge of the actual target domain layout is then discarded and, starting from a circle, we move the interior boundary $\partial\Omega_0$ with the intention to recover $z$ from the newly calculated $u$. 

The Lagrangian of this problem for $u\in u_0 + H_{(\Gamma_0\cup\Gamma_2)}^1(\Omega)$ and $\lambda\in H_{(\Gamma_0\cup\Gamma_2)}^1(\Omega)$ is 
\begin{align*}
  \Lag(u, \Omega, \lambda) = \frac{1}{2}\int_\Omega(u-z)^2\d x + \frac{\alpha}{2}\int_{\Omega_0}1 \d x + \mu\int_\Omega\nabla u^\top\nabla \lambda - f\lambda\d x.
\end{align*}
We use the principles presented in section \ref{mat-der-test} in order to derive first order shape derivatives. For ease of presetation, we choose the volume of $\Omega_0$ as regularization. Step 1 is already accomplished by building the Lagrangian above. This is now differentiated with respect to the shape.
\begin{align*}
  &\d{\Lag(u,\Omega,\lambda)}[\dm{u}{[\bV]}, \bV, \dm{\lambda}{[\bV]}] \\
  =& \int\limits_{\Omega} \div(\bV)\left(\frac{1}{2}(u-z)^2 
  + \frac{\alpha}{2}
  + \mu\inner{\nabla u}{\nabla \lambda} - \lambda f\right)
  \ \d \bx \\
   &+ \int\limits_\Omega \dm{\left(\frac{1}{2}(u-z)^2
   + \frac{\alpha}{2}
   + \mu\inner{\nabla u}{\nabla \lambda} - \lambda f\right)}{[\bV]} \ \d \bx \\
     =& \int\limits_{\Omega} \div(\bV)\left(\frac{1}{2}(u-z)^2
     + \frac{\alpha}{2}
     + \mu\inner{\nabla u}{\nabla \lambda} - \lambda f\right) \ \d \bx \\
   &+ \int\limits_\Omega (u-z)(\dm{u}{[\bV]} - \dm{z}{[\bV]}) + \mu\dm{\inner{\nabla u}{\nabla \lambda}}{[\bV]} + \inner{\nabla u}{\nabla \lambda}\dm{\mu}{[\bV]} \ \d \bx \\
   &+ \int\limits_{\Omega} - f\dm{\lambda}{[\bV]} - \lambda \dm{f}{[\bV]} \ \d \bx.
\end{align*}

Now, we use the well known identity
\begin{align*}
  \dm{\inner{\nabla u}{\nabla \lambda}}{[\bV]} = \inner{\nabla \dm{u}{[\bV]}}{\nabla\lambda} + \inner{\nabla u}{\nabla \dm{\lambda}{[\bV]}} - \inner{\nabla u}{(\D\bV + \D\bV^T)\nabla \lambda}.
\end{align*}
Thus, the Lagrangian can be rephrased as
\begin{equation}\label{eq:MixedFiniteElementLagrangianDerivative}
\begin{aligned}
&\d{\Lag(u,\Omega,\lambda)}{[\dm{u}{[\bV]}, \bV, \dm{\lambda}{[\bV]}]} \\
     =& \int\limits_{\Omega} \div(\bV)\left(\frac{1}{2}(u-z)^2
     + \frac{\alpha}{2}
     + \mu\inner{\nabla u}{\nabla \lambda} - \lambda f\right)
     - (u-z)\dm{z}{[\bV]} - \lambda \dm{f}{[\bV]}\\
     &\hspace{1cm} - \mu\inner{\nabla u}{(\D\bV + \D\bV^T)\nabla \lambda} + \inner{\nabla u}{\nabla \lambda}\dm{\mu}{[\bV]} \ \d \bx\\
   &+ \int\limits_\Omega (u-z)\dm{u}{[\bV]} + \mu\inner{\nabla \dm{u}{[\bV]}}{\nabla\lambda} + \mu\inner{\nabla u}{\nabla \dm{\lambda}{[\bV]}} - f\dm{\lambda}{[\bV]} \ \d \bx.
\end{aligned}
\end{equation}
The condition 
$\d{\Lag(u,\Omega,\lambda)}{[\dm{u}{[\bV]}, \bV, \dm{\lambda}{[\bV]}]} =0$ for all $\dm{\lambda}{[\bV]}\in H_{(\Gamma_0\cup\Gamma_2)}^1(\Omega)$ gives the state equation back. The 
condition $\dm{\Lag(u,\Omega,\lambda)}{[\dm{u}{[\bV]}, \bV, \dm{\lambda}{[\bV]}]} =0$ for all $\dm{u}{[\bV]}\in H_{(\Gamma_0\cup\Gamma_2)}^1(\Omega)$ results in the adjoint equation
\begin{align}\label{eq:MotherShapeAdjoint}
  \int\limits_\Omega \mu\inner{\nabla \dm{u}{[\bV]}}{\nabla\lambda} + (u-z)\dm{u}{[\bV]} \ \d \bx = 0 \quad \forall \dm{u}{[\bV]} \in H_{(\Gamma_0\cup\Gamma_2)}^1(\Omega).
\end{align}
By using the state and adjoint equation, we obtain the simplified shape derivative of the Lagrangian
\begin{equation}\label{foundations_dML}
\begin{aligned}
&\d{\Lag(u,\Omega,\lambda)}{[\bV]}\\
=&\int\limits_{\Omega} \div(\bV)\left(\frac{1}{2}(u-z)^2
 + \frac{\alpha}{2}
+ \mu\inner{\nabla u}{\nabla \lambda} - \lambda f 
\right) - (u-z)\dm{z}{[\bV]} - \lambda \dm{f}{[\bV]}\\
     &\hspace{1cm} - \mu\inner{\nabla u}{(\D\bV + \D\bV^T)\nabla \lambda} + \inner{\nabla u}{\nabla \lambda}\dm{\mu}{[\bV]} \ \d \bx.
\end{aligned}
\end{equation}
In the particular case of the numerical example, there is $f\equiv 0$, $z$ is a fixed field and $\mu$ deformes with the shape, and thus
$\dm{f}{[\bV]}=0, \dm{\mu}{[\bV]}=0, \dm{z}{[\bV]}=\nabla z^\top \bV$, which simplifies the shape derivative even more to
\begin{equation}\label{foundations_dML_sim}
\begin{aligned}
&\d{\Lag(u,\Omega,\lambda)}{[\bV]}\\
=&\int\limits_{\Omega} \div(\bV)\left(\frac{1}{2}(u-z)^2
 + \frac{\alpha}{2}
+ \mu\inner{\nabla u}{\nabla \lambda} - \lambda f 
\right) - (u-z)\nabla z^\top \bV  \\
     &\hspace{1cm} - \mu\inner{\nabla u}{(\D\bV + \D\bV^T)\nabla \lambda}  \ \d \bx.
\end{aligned}
\end{equation}
Here, we assume that $u$ and $\lambda$ satisfy the primal and adjoint equation.

The linear second shape derivative $J^{\prime\prime}(\Omega)[V,W]$ is derived in the appendix.

\section{Numerical Results}\label{sec:numres}
The Taylor series expansion shown in \Cref{algo_shape_Taylor} gives the potential of quadratic convergence for shape optimization algorithms using the linear second shape derivative introduced there, since the discussion before this theorem shows that standard shape optimization can be view from a linear vector space perspective. The performance obstacle is however the lack of positive definiteness of the linear second shape derivative. \Cref{remark-Tichonov-KKT} specifies a convenient strategy to circumvent this problem at the cost of loosing quadratic convergence, but getting almost arbitrarily good linear convergence. The aim of this section is to illustrate this effect in numerical computations for the model problem introduced in 
\Cref{sec:model}. Here, we first give more details on the algorithmic realization and afterwards illustrate the interplay of  \Cref{algo_shape_Taylor} and \Cref{remark-Tichonov-KKT}.

The variational Newton method aims at solving the stationarity condition for the Lagrangian with respect to all variables. Thus, it is a method iterating over all variables simultaneously in the following form 
\begin{align*}
  (u,\Omega,\lambda)_{k+1} = (u,\Omega,\lambda)_k + (\hat{u}, \hat{V}(\Omega_k), \hat{\lambda})\, ,
\end{align*}
where $\Omega_k+\hat{V}(\Omega_k)$ has to be read as $(I+\hat{V})(\Omega_k)$ and $(\hat{u}, \hat{V}, \hat{\lambda})$ solve the variational problem for the Lagrangian $\Lag$ defined in the previous section, discretized by continuous finite elements
\begin{align} \label{one-shot-var-KKT}
  \Lag^{\prime\prime}(u_{k},\Omega_{k},\lambda_{k})[\tilde{u}, \tilde{V}, \tilde{\lambda}][\hat{u}, \hat{V}, \hat{\lambda}] = &-\d \Lag(u_{k},\Omega_{k},\lambda_{k})[\tilde{u}, \tilde{V}, \tilde{\lambda}] \\ \nonumber
  &\forall (\tilde{u}, \tilde{q}, \tilde{\lambda}) \in \CG_{r_1} \times \CG^2_{r_2} \times \CG_{r_1},
\end{align}
Here, $\Lag^{\prime\prime}$ means the linear second order derivatives as in \cref{algo_shape_Taylor}, which is only of importance for the shape part, of course. In this way nonsymmetric terms of type $\tilde{V}\circ\hat{V}$ do not arise. As discussed above, this is just a particular formulation of a Newton method iterating over the state/adjoint variables and the deformation vector field of the shape. Thus, classical Newton convergence theory applies in linear spaces. 

Equation (\ref{one-shot-var-KKT}) is formulated for a particular and typical choice of ansatz and test spaces, which can be, of course, adapted to a specific application. Any solution algorithm has to cope with the fact the the shape Hessian has a huge kernel due to the Hadamard structure theorem. As already discussed, the simplest approach is to add a Tikhonov type regularization term involving a coercive bilinear form $b_{\Omega_k}$ such that the regularized variant

\begin{align} \label{one-shot-var-KKT-reg}
  \Lag^{\prime\prime}(u_{k},\Omega_{k},\lambda_{k})[\tilde{u}, \tilde{V}, \tilde{\lambda}][\hat{u}, \hat{V}, \hat{\lambda}] 
  + \frac{\varepsilon}{2}b_{\Omega_k}(\tilde{V},\hat{V})= &-\d \Lag(u_{k},\Omega_{k},\lambda_{k})[\tilde{u}, \tilde{V}, \tilde{\lambda}] \\ \nonumber
  &\forall (\tilde{u}, \tilde{q}, \tilde{\lambda}) \in \CG_{r_1} \times \CG^2_{r_2} \times \CG_{r_1},
\end{align}
possesses a unique solution, which converges to the Moore-Penrose pseudoinverse solution for $\varepsilon\to 0$, cf.~Remark \ref{remark-Tichonov-KKT}.

With obvious abbreviations, in particular $\Lag_\Omega$ and $\Lag_{\Omega\Omega}$ for first and linear second shape derivative, the variational KKT system (\ref{one-shot-var-KKT-reg}) can be written in matrix form as
\[
\mat[3]{
\Lag_{uu}&\Lag_{u\Omega} &\Lag_{u\lambda}\\
\Lag_{\Omega u}& \Lag_{\Omega\Omega}+\varepsilon b_{\Omega\Omega}&\Lag_{\Omega\lambda}\\
\Lag_{\lambda u}& \Lag_{\lambda\Omega}& 0
}
\mat(1){\hat{u}\\ \hat{V} \\\hat{\lambda}}=
-\mat(1){\Lag_u \\ \Lag_{\Omega} \\ \Lag_\lambda}
\]

%
This is the method, whose convergence history is shown the figures below after iteration 20. The first iterations are carried our by the following variant.
At the beginning of the nonlinear iterations, the second shape derivative may lack positive semidefiniteness. Then it is advisable to start with a method, which replace the Hessian terms by zero in the upper left quadrant, resulting in the equation

\[
\mat[3]{
0&0 &\Lag_{u\lambda}\\
0& \varepsilon b_{\Omega\Omega}&\Lag_{\Omega\lambda}\\
\Lag_{\lambda u}& \Lag_{\lambda\Omega}& 0
}
\mat(1){\hat{u}\\ \hat{V} \\\hat{\lambda}}=
-\mat(1){\Lag_u \\ \Lag_{\Omega} \\ \Lag_\lambda}
\]
which results in a projected steepest descent method  in the sense of \cite{Iliapap, Ilia-thesis}, which can conceptually also be enhanced by additional equality and inequality constraints. The latter can be reformulated variationally in the form

\begin{align} \label{one-shot-var-KKT-red}
  \Lag_{\lambda u}[\tilde{u}][\hat{\lambda}] + &\Lag_{u\lambda}[\tilde{\lambda}][\hat{u}] + \Lag_{\lambda\Omega}[\tilde{\lambda}][\hat{V}]
  +\Lag_{\Omega\lambda}[\tilde{V}][\hat{\lambda}]\\ \nonumber
  &+ \frac{\varepsilon}{2}b_{\Omega_k}(\tilde{V},\hat{V})= -\d \Lag(u_{k},\Omega_{k},\lambda_{k},\mu_k)[\tilde{u}, \tilde{V}, \tilde{\lambda}] \\ \nonumber
  &\forall (\tilde{u}, \tilde{q}, \tilde{\lambda}) \in \CG_{r_1} \times \CG^2_{r_2} \times \CG_{r_1},
\end{align}
thus giving rise to a convenient treatment of additional constraints also in the steepest descent case. This method is carried out in the first 5 iterations, in order to bring the method in the convergence vicinity of the shape Newton method. The regularization term is chosen in the following form
\[
b_{\Omega_k}({\bW},{\bV}) = \int\limits_\Omega \epsilon_1\left(\inner{\bW}{\bV} + \epsilon_2\inner{\nabla \bW}{\nabla \bV}_F\right) \ \d \bx
\]
with $\epsilon_2=5\cdot 10^{-1}$ in all our examples. The influence of the value of $\epsilon_1$ is the main discussion aspect below.

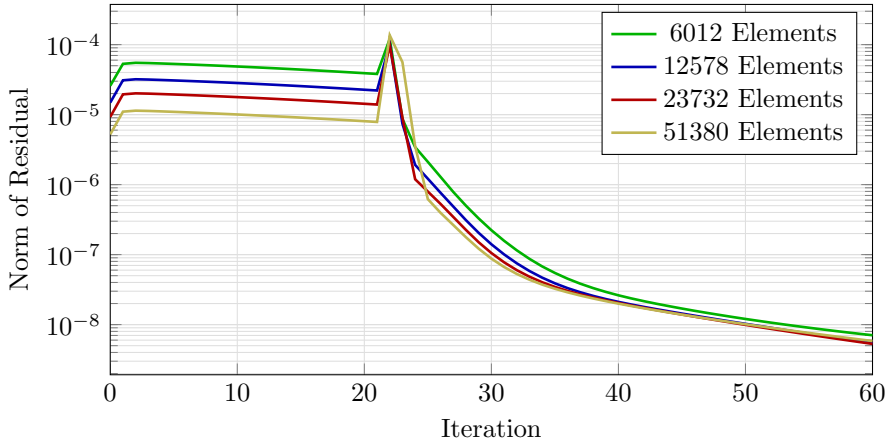
\begin{figure}[h!]
\begin{center}
\begin{tikzpicture}
\begin{axis}[
	ymode=log,
	xmin=0.0,
	xmax=60.0,
	xlabel=Iteration, 
	ylabel=Norm of Residual,
	grid=both,
	minor grid style={gray!25},
	major grid style={gray!25},
	width=0.9\linewidth,
	height=0.5\linewidth,
	no marks]
\addplot[line width=1pt,solid,color=green!70!black] table[x index = {0}, y index = {3},col sep=space]{history_circle_1_9.txt};\addlegendentry{$6012$ Elements};
\addplot[line width=1pt,solid,color=blue!70!black] table[x index = {0}, y index = {3},col sep=space]{history_circle_2_9.txt};\addlegendentry{$12578$ Elements};
\addplot[line width=1pt,solid,color=red!70!black] table[x index = {0}, y index = {3},col sep=space]{history_circle_3_9.txt};\addlegendentry{$23732$ Elements};
\addplot[line width=1pt,solid,color=yellow!70!black] table[x index = {0}, y index = {3},col sep=space]{history_circle_4_9.txt};\addlegendentry{$51380$ Elements};
\end{axis}
\end{tikzpicture}
\end{center}
\caption{Effects of the Nowton-type algorithm on different meshes. Here, always $\epsilon_1=3\cdot 10^{-2}$.}\label{fig:MeshConvergence}
\end{figure}

\Cref{fig:MeshConvergence} shows the convergence history of the Newton-like iterations for $\epsilon_1=3\cdot 10^{-2}$ on different meshes. The first 20 iterations are done by the projected gradient strategy formulated at the beginning of this section. In these first 20 iterations, the step size was chosen \STS{$0.9$}. Starting from iteration 20, in the Newton part, only full steps with stepsize 1 have been performed. After iteration 20, we observe a very fast convergence, which levels off later to still good linear convergence. 

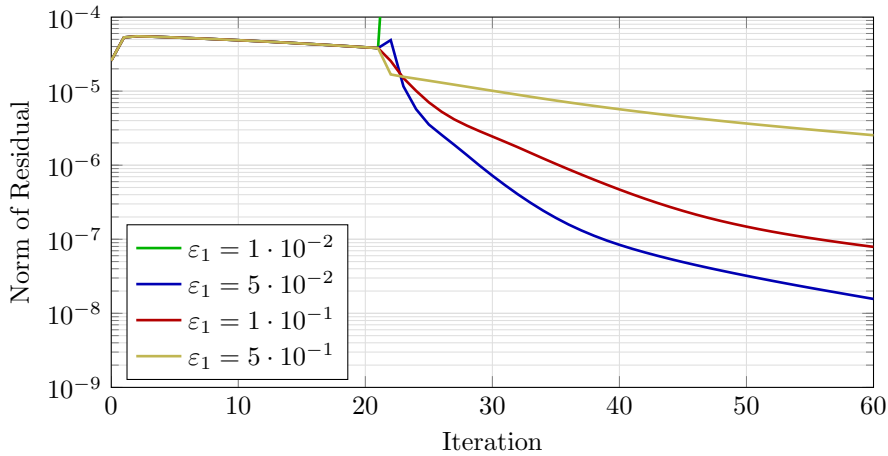
\begin{figure}[h!]
\begin{center}
\begin{tikzpicture}
\begin{axis}[
	ymode=log,
	xmin=0.0,
	xmax=60.0,
	ymin=1e-9,
	ymax=1e-4,
	xlabel=Iteration, 
	ylabel=Norm of Residual,
	grid=both,
	minor grid style={gray!25},
	major grid style={gray!25},
	width=0.9\linewidth,
	height=0.5\linewidth,
	legend style={at={(0.02,0.02)},anchor=south west},
	no marks]
\addplot[line width=1pt,solid,color=green!70!black] table[x index = {0}, y index = {3},col sep=space]{history_m20_01.txt};\addlegendentry{$\varepsilon_1 = 1\cdot10^{-2}$};
\addplot[line width=1pt,solid,color=blue!70!black] table[x index = {0}, y index = {3},col sep=space]{history_m20_05.txt};\addlegendentry{$\varepsilon_1 = 5\cdot 10^{-2}$};
\addplot[line width=1pt,solid,color=red!70!black] table[x index = {0}, y index = {3},col sep=space]{history_m20_1.txt};\addlegendentry{$\varepsilon_1 = 1\cdot10^{-1}$};
\addplot[line width=1pt,solid,color=yellow!70!black] table[x index = {0}, y index = {3},col sep=space]{history_m20_5.txt};\addlegendentry{$\varepsilon_1 = 5\cdot10^{-1}$};
\end{axis}
\end{tikzpicture}
\end{center}
\caption{Convergence history for mesh study for the choise of different values for $\epsilon_1$ on a mesh with 6012 elements.}\label{fig:eps-1-study-1}
\end{figure}

\begin{figure}[h!]
\begin{center}
\begin{tikzpicture}
\begin{axis}[
	ymode=log,
	xmin=0.0,
	xmax=60.0,
	ymin=1e-9,
	ymax=1e-4,
	xlabel=Iteration, 
	ylabel=Norm of Residual,
	grid=both,
	minor grid style={gray!25},
	major grid style={gray!25},
	width=0.9\linewidth,
	height=0.5\linewidth,
	legend style={at={(0.02,0.02)},anchor=south west},
	no marks]
\addplot[line width=1pt,solid,color=green!70!black] table[x index = {0}, y index = {3},col sep=space]{history_m20_4_9_01.txt};\addlegendentry{$\varepsilon_1 = 1\cdot10^{-2}$};
\addplot[line width=1pt,solid,color=blue!70!black] table[x index = {0}, y index = {3},col sep=space]{history_m20_4_9_05.txt};\addlegendentry{$\varepsilon_1 = 5\cdot 10^{-2}$};
\addplot[line width=1pt,solid,color=red!70!black] table[x index = {0}, y index = {3},col sep=space]{history_m20_4_9_1.txt};\addlegendentry{$\varepsilon_1 = 1\cdot10^{-1}$};
\addplot[line width=1pt,solid,color=yellow!70!black] table[x index = {0}, y index = {3},col sep=space]{history_m20_4_9_5.txt};\addlegendentry{$\varepsilon_1 = 5\cdot10^{-1}$};
\end{axis}
\end{tikzpicture}
\end{center}
\caption{Convergence history for mesh study for the choise of different values for $\epsilon_1$ on a mesh with $51380$ elements.}\label{fig:eps-1-study-2}
\end{figure}
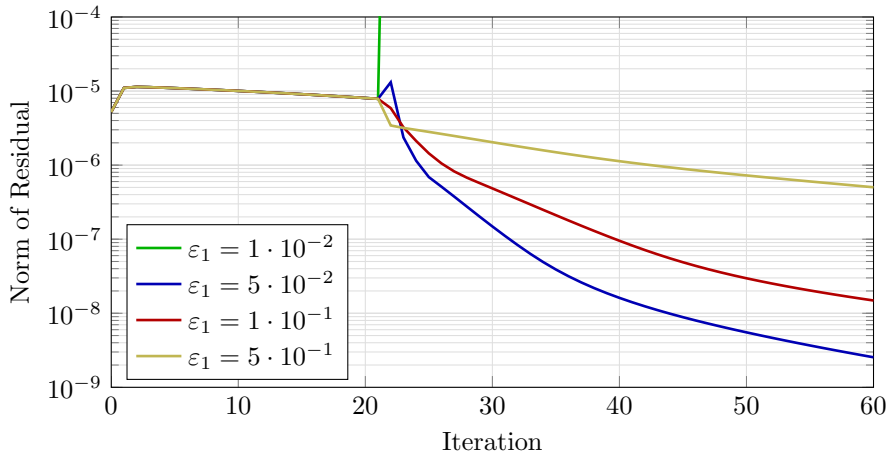

\Cref{fig:eps-1-study-1} and \Cref{fig:eps-1-study-2} investigate the effect of the particular choice of $\epsilon_1$ on the overall convergence on two different meshes. We observe that the convergence in the Newton part of the algorithm can be dramatically improved by reducing the parameter $\epsilon_1$, as is to be expected by the previous discussions in the paper. However, there is, of course, a lower limit for this parameter, which has to be larger than zero. In both figures, we observe divergence of the algorithm for $\epsilon = 10^{-2}$ and very good convergence for the slightly higher value $\epsilon=5\cdot 10^{-2}$.

All computations were performed on unstructured meshes, where the shape optimization algorithms have been implemented by using the finite elements toolbox  FEniCS \cite{AlnaesBlechta2015a}.

\section{Conclusions}
\label{sec:conclusions}
This papers analyses standard descending shape optimization algorithms from the point of view of iterating over deformations rather than geometries. It turns out that standard choices for scalar products in Hilbert spaces of deformations show indeed equivalence of shape optimization algorithms to optimization in a Hilbert space of deformations. Based on that, a linear variant of the second shape derivative is derived, on which one can base a Taylor expansion in linear spaces as well as a standard Newton iteration. The central obstacle of this point of view is the significant rank deficiency of the linear second shape derivative. As a remedy, Newton-type methods based on pseudoinverses are presented, as well as,  a convenient way to substitute the otherwise numerically challenging pseudoinverse. 

Finally, the question remains, whether we should replace standard shape optimization algorithms in the form of \cref{algo-shape-steepest} by a rather more straight forward iteration in the linear Hilbert space of deformations. The catch is that all intermediate and also the final domain deformations have to be invertible. This property can be much easier checked in each iteration of  \cref{algo-shape-steepest} separately rather than employing a standard formulation of a descent algorithm in vector spaces. Thus, we still recommend algorithms of the type \cref{algo-shape-steepest}, but to keep  in mind that the linear point of view provides a very convenient tool for the analysis of these algorithms.

\section*{Acknowledgments}
Both authors acknowledge support by the Deutsche Forsch\-ungsgemeinschaft within the Priority program SPP 1962 ``Non-smooth and Complem\-entarity-based Distributed Parameter Systems: Simulation and Hierarchical Optimization'' and the second author also within the Research Training Group RTG 2126 ``Algorithmic Optimization''.

\appendix
\section{Second Derivatives for the Interface Problem}
Here, the full Newton KKT system is derived in variational form. From \cref{sec:main} it is clear that we need the linear second shape derivative $\Lag^{\prime\prime}$. This is, however derived in the following manner: we investigate the standard second shape derivative $d^2\Lag$ and exclude the nonsymmetric terms on the way.
Indeed, in order to derive the KKT-system, the steps for the first shape derivative can very much be repeated to derive higher order shape derivatives. The volume formulation thus makes considering the Hessian almost as elegant as deriving the gradient.

Recall that in a general setting, the volume form of the shape derivative for a volume objective is given by
\begin{align*}
\d{J(\Omega)}{[\bV]} = \int\limits_{\Omega} g\div\bV + \dm{g}{[\bV]} \ \d \bx.
\end{align*}
Applying the same transformation again, one arrives at
\begin{align}\label{eq:VolHess1}
  \d^2{J}{[\bV, \bW]} = \int\limits_\Omega \left(g\div\bV + \dm{g}{[\bV]}\right)\div\bW + \dm{\left(g\div\bV + \dm{g}{[\bV]}\right)}{[\bW]} \ \d \bx.
\end{align}
Using the rule $\dm{(\D u)}{[\bV]} =\D(\dm{u}{[\bV]}) - \D u\D\bV$ for any function $u$, one can see that
\begin{align*}
  \dm{\left(\div\bV\right)}{[\bW]} &= \dm{\left(\tr\D\bV\right)}{[\bW]} = \tr\left(\dm{\left(\D\bV\right)}{[\bW]}\right) = \tr\left(\D\dm{\bV}{[\bW]} - \D\bV\D\bW\right)\\
  &= \div\left(\dm{\bV}{[\bW]}\right) - \tr(\D\bV\D\bW).
\end{align*}
Within the finite element context, one can consider the mesh motions $\bV$ and $\bW$ to be given as Vector-Lagrange functions in $\CG^3(\mathcal{D})$. Hence, they are automatically transported with the mesh and one immediately arrives at $\dm{\bV}{[\bW]} = 0$. This simplifies the above within the finite element context to
\begin{align*}
  \dm{\left(\div\bV\right)}{[\bW]} &= - \tr(\D\bV\D\bW).
\end{align*}
Inserting everything into~\eqref{eq:VolHess1} and rearranging leads to the very compact and elegant representation of the shape Hessian for a volume integral:
\begin{equation}\label{eq:GeneralVolumeHessian}
\begin{aligned}
  \d^2{J}{[\bV, \bW]} = &\int\limits_\Omega g\cdot\left(\div(\bV)\div(\bW) - \tr\left(\D\bV\D\bW\right)\right) \\
  &\hspace{1cm} + \dm{g}{[\bV]}\div(\bW) + \dm{g}{[\bW]}\div(\bV) + \ddm{g}{[\bV,\bW]} \ \d \bx.
\end{aligned}
\end{equation}

The shape derivative procedure can now be applied to the material derivative of the Lagrangian in~\eqref{eq:MixedFiniteElementLagrangianDerivative}, which will then lead to the KKT-system in variational form, ready to be implemented using mixed finite elements.
\begin{equation}\label{eq:LagrangianMotherRecall}
\begin{aligned}
&\d{\Lag(u,\Omega,\lambda)}{[\dm{u}{[\bV]}, \bV, \dm{\lambda}{[\bV], \tilde{\lambda}_{Vol}}]} \\
     =& \int\limits_{\Omega} \div(\bV)\left(\frac{1}{2}(u-z)^2
      + \frac{\alpha}{2}
     + \mu\inner{\nabla u}{\nabla \lambda} - \lambda f + \lambda_{Vol}\right)\\
   &\hspace{1cm} - (u-z)\dm{z}{[\bV]} - \lambda \dm{f}{[\bV]}\\
   &\hspace{1cm} - \mu\inner{\nabla u}{(\D\bV + \D\bV^T)\nabla \lambda} + \inner{\nabla u}{\nabla \lambda}\dm{\mu}{[\bV]} \\
   &\hspace{1cm}+ (u-z)\dm{u}{[\bV]} + \inner{\nabla \dm{u}{[\bV]}}{\nabla\lambda} + \inner{\nabla u}{\nabla \dm{\lambda}{[\bV]}} - f\dm{\lambda}{[\bV]} \ \d \bx
\end{aligned}
\end{equation}
Applying~\eqref{eq:GeneralVolumeHessian} to~\eqref{eq:LagrangianMotherRecall} will thus create an excessive amount of terms. To keep the derivation somewhat readable, we split~\eqref{eq:GeneralVolumeHessian} into separate terms
\begin{align*}
  T_1[\bV, \bW] \coloneq& g\cdot\left(\div(\bV)\div(\bW) - \tr\left(\D\bV\D\bW\right)\right)\\
  T_2[\bV, \bW] \coloneq& \dm{g}{[\bV]}\div(\bW)\\
  T_3[\bV, \bW] \coloneq& \ddm{g}{[\bV,\bW]}
\end{align*}
The KKT-system can then be written as
\begin{align*}
  \ddm{\Lag}{[\bV, \bW]} = \int\limits_{\Omega} T_1[\bV, \bW] + T_2[\bV, \bW] + T_2[\bW, \bV] + T_3[\bV, \bW] \ \d \bx.
\end{align*}
From the derivation of the gradient, we recall that in this special setting, we have
\begin{align*}
  g &= \frac{1}{2}(u-z)^2
   + \frac{\alpha}{2}
  + \mu\inner{\nabla u}{\nabla \lambda} - \lambda f
\end{align*}
and
\begin{align*}
\dm{g}{[\bV]} =& -(u-z)\dm{z}{[\bV]} - \lambda \dm{f}{[\bV]} - \mu\inner{\nabla u}{(\D\bV + \D\bV^T)\nabla \lambda} \\
   &+ (u-z)\dm{u}{[\bV]} + \mu\inner{\nabla \dm{u}{[\bV]}}{\nabla\lambda} + \mu\inner{\nabla u}{\nabla \dm{\lambda}{[\bV]}}\\
   &- f\dm{\lambda}{[\bV]}.
\end{align*}
We now have to formally derive $\ddm{g}{[\bV,\bW]} := \dm{\left(\dm{g}{[\bV]}\right)}{[\bW]}$, which again generates an excessive amount of terms. To keep everything readable, we split the above second material derivative into even smaller parts, namely
\begin{align*}
  T_{3,2} \coloneq& -(u-z)\dm{z}{[\bV]} - \lambda \dm{f}{[\bV]}\\
\Rightarrow \dm{T_{3,2}}{[\bW]} =& -(\dm{u}{[\bW]}-\dm{z}{[\bW]})\dm{z}{[\bV]} - (u-z)\ddm{z}{[\bV,\bW]} \\
  &-\dm{\lambda}{[\bW]}\dm{f}[\bV] - \lambda\ddm{f}{[\bV, \bW]}
\end{align*}
and
\begin{align*}
 T_{3,3} \coloneq &\mu\inner{\nabla u}{(\D\bV + \D\bV^T)\nabla \lambda}
\end{align*}
and the material derivative of this expression is computed as follows
\begin{align*}
  \dm{T_{3,3}}{[\bW]} =& - \mu\inner{\dm{\nabla u}{[\bW]}}{(\D\bV + \D\bV^T)\nabla \lambda} - \mu\inner{\nabla u}{\dm{\left((\D\bV + \D\bV^T)\nabla \lambda\right)}{[\bW]}} \\
  =& - \mu\inner{\nabla \dm{u}{[\bW]} - \D\bW^T \nabla u}{(\D\bV + \D\bV^T)\nabla \lambda} \\
  &- \mu\inner{\nabla u}{\dm{\left(\D\bV + \D\bV^T\right)}{[\bW]}\nabla \lambda} - \mu\inner{\nabla u}{(\D\bV + \D\bV^T)\dm{(\nabla \lambda)}{[\bW]}}\\
  =& -\mu\inner{\nabla \dm{u}{[\bW]}}{(\D\bV + \D\bV^T)\nabla \lambda} + \mu\inner{\D\bW^T \nabla u}{(\D\bV + \D\bV^T)\nabla \lambda} \\
  &- \mu\inner{\nabla u}{(\D\dm{\bV}{[\bW]} - \D\bV\D\bW + \D\dm{\bV}{[\bW]}^T - (\D\bV\D\bW)^T)\nabla \lambda} \\
  &- \mu\inner{\nabla u}{(\D\bV + \D\bV^T)\nabla\dm{\lambda}{[\bW]}} + \mu\inner{\nabla u}{(\D\bV + \D\bV^T)\D\bW^T \nabla\lambda}\\
    =& -\mu\inner{\nabla \dm{u}{[\bW]}}{(\D\bV + \D\bV^T)\nabla \lambda} + \mu\inner{\D\bW^T \nabla u}{(\D\bV + \D\bV^T)\nabla \lambda} \\
  &+ \mu\inner{\nabla u}{(\D\bV\D\bW + (\D\bV\D\bW)^T)\nabla \lambda} \\
  &- \mu\inner{\nabla u}{(\D\bV + \D\bV^T)\nabla\dm{\lambda}{[\bW]}} + \mu\inner{\nabla u}{(\D\bV + \D\bV^T)\D\bW^T \nabla\lambda},
\end{align*}
where we have again used $\dm{\bV}[\bW] = 0$ in the last step. This object can be simplified further
\begin{align*}
  \dm{T_{3,3}}{[\bW]} =& -\mu\inner{\nabla \dm{u}{[\bW]}}{(\D\bV + \D\bV^T)\nabla \lambda} - \mu\inner{\nabla u}{(\D\bV + \D\bV^T)\nabla\dm{\lambda}{[\bW]}}\\
  &+ \mu\inner{\nabla u}{\D\bW(\D\bV + \D\bV^T)\nabla \lambda} + \mu\inner{\nabla u}{(\D\bV + \D\bV^T)\D\bW^T \nabla\lambda}\\
  &+ \mu\inner{\nabla u}{(\D\bV\D\bW + (\D\bV\D\bW)^T)\nabla \lambda} \\
  =& -\mu\inner{\nabla \dm{u}{[\bW]}}{(\D\bV + \D\bV^T)\nabla \lambda} - \mu\inner{\nabla u}{(\D\bV + \D\bV^T)\nabla\dm{\lambda}{[\bW]}}\\
  &+ \mu\inner{\nabla u}{(\D\bW\D\bV + \D\bW\D\bV^T)\nabla \lambda} + \mu\inner{\nabla u}{(\D\bV\D\bW^T + (\D\bW\D\bV)^T) \nabla\lambda}\\
  & + \mu\inner{\nabla u}{(\D\bV\D\bW + (\D\bV\D\bW)^T)\nabla \lambda}\\
  =& -\mu\inner{\nabla \dm{u}{[\bW]}}{(\D\bV + \D\bV^T)\nabla \lambda} - \mu\inner{\nabla u}{(\D\bV + \D\bV^T)\nabla\dm{\lambda}{[\bW]}}\\
  &+ \mu\inner{\nabla u}{(\D\bW\D\bV + \D\bV\D\bW)\nabla\lambda}\\
  &+ \mu\inner{\nabla u}{(\D\bW\D\bV^T + \D\bV\D\bW^T)\nabla \lambda}\\
  &+ \mu\inner{\nabla u}{(\D\bW\D\bV)^T + (\D\bV\D\bW)^T) \nabla\lambda}.
\end{align*}

The material derivative of the next component is again straight forward
\begin{align*}
  T_{3,4} \coloneq& (u-z)\dm{u}{[\bV]}\\
 \dm{T_{3,4}}{[\bW]} =&\dm{\left((u-z)\dm{u}{[\bV]}\right)}{[\bW]}\\
  = &(\dm{u}{[\bW]} - \dm{z}{[\bW]})\dm{u}{[\bV]} + (u-z)\dm{\left(\dm{u}{[\bV]}\right)}{[\bW]}.
\end{align*}
The term $\dm{\left(\dm{u}{[\bV]}\right)}{[\bW]}$ belongs to the purely nonsymmetric part the second shape derivative and is left out in the linear second shape derivative. Thus, instead of $\dm{T_{3,4}}{[\bW]}$, we use later on 
\[ 
\dm{T_{3,4}^\text{lin}}{[\bW]}  = (\dm{u}{[\bW]} - \dm{z}{[\bW]})\dm{u}{[\bV]}.
\]
The second to last expression is given by
\begin{align*}
T_{3,5} \coloneq& \mu\inner{\nabla \dm{u}{[\bV]}}{\nabla\lambda} + \mu\inner{\nabla u}{\nabla \dm{\lambda}{[\bV]}}
\end{align*}
and the material derivative is computed as follows
\begin{align*}
\dm{T_{3,5}}{[\bW]} =&\mu\inner{\dm{\left(\nabla \dm{u}{[\bV]}\right)}{[\bW]}}{\nabla \lambda} + \mu\inner{\nabla \dm{u}{[\bV]}}{\dm{\left(\nabla \lambda\right)}{[\bW]}} \\
&+ \mu\inner{\dm{\left(\nabla u\right)}{[\bW]}}{\nabla \dm{\lambda}{[\bV]}} + \mu\inner{\nabla u}{\dm{\left(\nabla\dm{\lambda}{[\bV]}\right)}{[\bW]}}\\
=& \mu\inner{\nabla\left(\dm{\left(\dm{u}{[\bV]}\right)}{[\bW]}\right)-\D \bW^T\nabla \dm{u}{[\bV]}}{\nabla \lambda}\\
&+ \mu\inner{\nabla \dm{u}{[\bV]}}{\nabla\left(\dm{\lambda}{[\bW]}\right) - \D\bW^T\nabla\lambda} \\
&+ \mu\inner{\nabla\left(\dm{u}{[\bW]}\right) - \D\bW^T \nabla u}{\nabla \dm{\lambda}{[\bV]}}\\
&+ \mu\inner{\nabla u}{\nabla\left(\dm{\left(\dm{\lambda}{[\bV]}\right)}{[\bW]}\right) - \D\bW^T\nabla\dm{\lambda}{[\bV]}}.
\end{align*}
Again, we ignore the nonsymmetric contributions $\dm{\left(\dm{u}{[\bV]}\right)}{[\bW]}$ and $\dm{\left(\dm{\lambda}{[\bV]}\right)}{[\bW]}$ in the changed term
\begin{align*}
\dm{T_{3,5}^\text{lin}}{[\bW]}=& -\mu\inner{\D \bW^T\nabla \dm{u}{[\bV]}}{\nabla \lambda} + \mu\inner{\nabla \dm{u}{[\bV]}}{\nabla\left(\dm{\lambda}{[\bW]}\right) - \D\bW^T\nabla\lambda} \\
&+ \mu\inner{\nabla\left(\dm{u}{[\bW]}\right) - \D\bW^T \nabla u}{\nabla \dm{\lambda}{[\bV]}} - \mu\inner{\nabla u}{\D\bW^T\nabla\dm{\lambda}{[\bV]}}
\end{align*}
This expression can also be simplified
\begin{align*}
  \dm{T_{3,5}^\text{lin}}{[\bW]} =& -\mu\inner{\nabla \dm{u}{[\bV]}}{(\D \bW+\D\bW^T)\nabla \lambda} + \mu\inner{\nabla \dm{u}{[\bV]}}{\nabla\dm{\lambda}{[\bW]}} \\
&+ \mu\inner{\nabla\dm{u}{[\bW]}}{\nabla \dm{\lambda}{[\bV]}} - \mu\inner{\nabla u}{(\D\bW + \D\bW^T)\nabla \dm{\lambda}{[\bV]}}.
\end{align*}
Last, but not least, we have
\begin{align*}
  T_{3,6} \coloneq& - f\dm{\lambda}{[\bV]}\\
  \dm{T_{3,6}}{[\bW]} =& -\dm{f}{[\bW]}\dm{\lambda}{[\bV]} - f\dm{\left(\dm{\lambda}{[\bV]}\right)}{[\bW]} 
\end{align*}
and thus
\[
  		\dm{T_{3,6}^\text{lin}}{[\bW]}		= -\dm{f}{[\bW]}\dm{\lambda}{[\bV]}.
\]

Taking everything together,
we arrive at a variational expression for the shape-KKT system, ready to be implemented using mixed finite elements
\begin{align*}
  &\Lag^{\prime\prime}[\dm{u}{[\bV]}, \bV, \dm{\lambda}{[\bV]}][\dm{u}{[\bW]}, \bW, \dm{\lambda}{[\bW]]} \\
  =& \int\limits_{\Omega} \left(\frac{1}{2}(u-z)^2
   + \frac{\alpha}{2}
  + \mu\inner{\nabla u}{\nabla \lambda} - \lambda f + \lambda_{Vol}\right)\left(\div(\bV)\div(\bW) - \tr\left(\D\bV\D\bW\right)\right) \\
   & +\left(-(u-z)\dm{z}{[\bV]} - \lambda \dm{f}{[\bV]} - \mu\inner{\nabla u}{(\D\bV + \D\bV^T)\nabla \lambda}\right)\div\bW \\
   & +\left(-(u-z)\dm{z}{[\bW]} - \lambda \dm{f}{[\bW]} - \mu\inner{\nabla u}{(\D\bW + \D\bW^T)\nabla \lambda}\right)\div\bV \\
   & +\left((u-z)\dm{u}{[\bV]} + \mu\inner{\nabla \dm{u}{[\bV]}}{\nabla\lambda} + \mu\inner{\nabla u}{\nabla \dm{\lambda}{[\bV]}} - f\dm{\lambda}{[\bV]}\right)\div\bW\\
   &+\left((u-z)\dm{u}{[\bW]} + \mu\inner{\nabla \dm{u}{[\bW]}}{\nabla\lambda} + \mu\inner{\nabla u}{\nabla \dm{\lambda}{[\bW]}} - f\dm{\lambda}{[\bW]}\right)\div\bV\\
   &+\dm{z}{[\bW]}\dm{z}{[\bV]} - (u-z)\ddm{z}{[\bV,\bW]} \\
   &+\dm{u}{[\bW]}\dm{u}{[\bV]} - \dm{u}{[\bV]}\dm{z}{[\bW]} - \dm{u}{[\bW]}\dm{z}{[\bV]}\\
   &-\dm{\lambda}{[\bW]}\dm{f}[\bV] -\dm{\lambda}{[\bV]}\dm{f}{[\bW]} - \lambda\ddm{f}{[\bV, \bW]}\\
   & -\mu\inner{\nabla \dm{u}{[\bW]}}{(\D\bV + \D\bV^T)\nabla \lambda} - \mu\inner{\nabla u}{(\D\bV + \D\bV^T)\nabla\dm{\lambda}{[\bW]}}\\
   & -\mu\inner{\nabla \dm{u}{[\bV]}}{(\D \bW+\D\bW^T)\nabla \lambda} - \mu\inner{\nabla u}{(\D\bW + \D\bW^T)\nabla \dm{\lambda}{[\bV]}} \\
  &+ \mu\inner{\nabla u}{(\D\bW\D\bV + \D\bV\D\bW)\nabla\lambda}\\
  &+ \mu\inner{\nabla u}{(\D\bW\D\bV^T + \D\bV\D\bW^T)\nabla \lambda}\\
  &+ \mu\inner{\nabla u}{(\D\bW\D\bV)^T + (\D\bV\D\bW)^T) \nabla\lambda}\\
&+ \mu\inner{\nabla\dm{u}{[\bW]}}{\nabla \dm{\lambda}{[\bV]}} + \mu\inner{\nabla \dm{u}{[\bV]}}{\nabla\dm{\lambda}{[\bW]}} \ \d\bx
\end{align*}
The term $\ddm{z}{[\bV,\bW]}=V^\top \text{Hess}(z)W$ has to be assembled using finite elements. Often $z$ itself is just a finite element approximation. If the order of this approximation is not higher than 1, then the assembled term is zero. 
Furthermore $\ddm{f}{[\bV, \bW]}=0$ since the right hand side $f$ is assumed to deform with the mesh.

\bibliographystyle{siamplain}

\end{document}